# HARRIS RECURRENCE OF METROPOLIS-WITHIN-GIBBS AND TRANS-DIMENSIONAL MARKOV CHAINS

By Gareth O. Roberts and Jeffrey S. Rosenthal

*Lancaster University and University of Toronto*

A $\phi$-irreducible and aperiodic Markov chain with stationary probability distribution will converge to its stationary distribution from almost all starting points. The property of Harris recurrence allows us to replace "almost all" by "all," which is potentially important when running Markov chain Monte Carlo algorithms. Full-dimensional Metropolis–Hastings algorithms are known to be Harris recurrent. In this paper, we consider conditions under which Metropolis-within-Gibbs and trans-dimensional Markov chains are or are not Harris recurrent. We present a simple but natural two-dimensional counter-example showing how Harris recurrence can fail, and also a variety of positive results which guarantee Harris recurrence. We also present some open problems. We close with a discussion of the practical implications for MCMC algorithms.

**1. Introduction.** Harris recurrence is a concept introduced fifty years ago by Harris [8]. More recently, connections between Harris recurrence and Markov chain Monte Carlo (MCMC) algorithms were investigated by Tierney [25] and Chan and Geyer [3]. In this paper, we re-examine Harris recurrence of various MCMC algorithms in a more general context.

Markov chains with stationary distributions are the basis of MCMC algorithms. For the algorithm to be valid, it is crucial that the chain converges to stationarity in distribution. If the state space is countable and the Markov chain is aperiodic and also (*classically*) *irreducible* (i.e., has positive probability of reaching any state from any other state), then it is well known that convergence to stationarity is guaranteed from all starting states (see, e.g., [1, 10, 22]).

On the other hand, classical irreducibility is unachievable when the state space $\mathcal{X}$ is uncountable. A weaker property is $\phi$-*irreducibility* [i.e., having positive probability of reaching every subset $A$ with $\phi(A) > 0$ from every state $x \in \mathcal{X}$, for some nonzero measure $\phi(\cdot)$]. It is known that a $\phi$-irreducible, aperiodic Markov chain with stationary probability distribution $\pi(\cdot)$ must still converge to $\pi(\cdot)$ from $\pi$-almost every starting point (see, e.g., [13, 15, 23]). However, if a chain is $\phi$-irreducible but not classically irreducible, then it is indeed possible to have a null set of states from which convergence does not occur.

Tierney [25] and Chan and Geyer [3] note that this null set of points from which convergence fails could cause practical problems for MCMC algorithms if the user happens to choose an initial state in this null set. Thus, understanding the nature of this null set is important for applications of MCMC, as well as theoretically. Chan and Geyer [3] refer to this null set as a "measure-theoretic pathology." However, we shall see herein that the null set can arise quite naturally, on both discrete and continuous state spaces, including for a simple two-dimensional Metropolis-within-Gibbs algorithm with continuous densities.

This paper is structured as follows. Section 2 presents some background about Markov chains and Harris recurrence and Theorem 6 proves a number of equivalences of Harris recurrence. Section 3 discusses full-dimensional Metropolis–Hastings algorithms and Section 4 then discusses Metropolis-within-Gibbs algorithms. Example 9 demonstrates that a simple two-dimensional Metropolis-within-Gibbs algorithm with continuous target and proposal densities, although irreducible and aperiodic, can still fail to be Harris recurrent or to converge to stationarity from all starting points. Sections 4 and 5 prove various positive results which guarantee Harris recurrence for Metropolis-within-Gibbs algorithms under various conditions and Section 6 does the same for trans-dimensional Markov chains.

**2. Markov chains and Harris recurrence.** Consider a Markov chain $\{X_n\}$ with transition probabilities $P(x, \cdot)$, on a state space $\mathcal{X}$ with $\sigma$-algebra $\mathcal{F}$. Let $P^n(x, \cdot)$ be the $n$-step transition kernel and for $A \in \mathcal{F}$, let $\tau_A = \inf\{n \geq 1 : X_n \in A\}$ be the first return time to $A$, with $\tau_A = \infty$ if the chain never returns to $A$.

Recall that a Markov chain is $\phi$-*irreducible* if there exists a nonzero $\sigma$-finite measure $\psi(\cdot)$ on $(\mathcal{X}, \mathcal{F})$ such that $\mathbf{P}[\tau_A < \infty | X_0 = x] > 0$ for all $x \in \mathcal{X}$ and all $A \in \mathcal{F}$ with $\psi(A) > 0$. The probability distribution $\pi(\cdot)$ on $(\mathcal{X}, \mathcal{F})$ is *stationary* for the chain if $\int_\mathcal{X} \pi(dx) P(x, A) = \pi(A)$ for all $A \in \mathcal{F}$.

The *period* of a $\phi$-irreducible chain with stationary distribution $\pi(\cdot)$ is the largest $D \in \mathbf{N}$ (the set of all positive integers) for which there exist disjoint subsets $\mathcal{X}_1, \mathcal{X}_2, \ldots, \mathcal{X}_D \in \mathcal{F}$ with $\pi(\mathcal{X}_i) > 0$, such that $P(x, \mathcal{X}_{i+1}) = 1$ for all $x \in \mathcal{X}_i$ ($1 \leq i \leq D-1$) and $P(x, \mathcal{X}_1) = 1$ for all $x \in \mathcal{X}_D$. If $D = 1$, then the chain is *aperiodic*.

In terms of these definitions, we have the following classical result, as in, for example, [25], page 1758 or [23]. (See also [15] and [13].)

PROPOSITION 1. *Consider a $\phi$-irreducible, aperiodic Markov chain with stationary probability distribution $\pi(\cdot)$. Let $G$ be the set of $x \in \mathcal{X}$ such that $\lim_{n\to\infty} \|P^n(x,\cdot) - \pi(\cdot)\| = 0$. Then $\pi(G) = 1$.*

We also note that aperiodicity is not essential in Proposition 1:

PROPOSITION 2. *Consider a $\phi$-irreducible Markov chain with stationary probability distribution $\pi(\cdot)$ and period $D \geq 1$. Let $G$ be the set of $x \in \mathcal{X}$ such that $\lim_{n\to\infty} \|(1/D)\sum_{r=1}^{D} P^{nD+r}(x,\cdot) - \pi(\cdot)\| = 0$. Then $\pi(G) = 1$.*

PROOF. If $D = 1$, then this reduces to Proposition 1 above. If $D > 1$, then let $\mathcal{X}_1, \mathcal{X}_2, \ldots, \mathcal{X}_D$ be as in the above definition of period. Then $P^D$ is $\phi$-irreducible and aperiodic when restricted to each $\mathcal{X}_i$, with stationary distribution $\pi_i(\cdot)$ such that $\pi(\cdot) = (1/D)\sum_{r=1}^{D} \pi_i(\cdot)$. It follows from Proposition 1 that for $\pi_D$-a.e. $x \in \mathcal{X}_D$ and for any $1 \leq r \leq D$, we have $\lim_{n\to\infty} \|P^{nD+r}(x,\cdot) - \pi_r(\cdot)\| = 0$. Hence, for $\pi_i$-a.e. $x \in \mathcal{X}_i$, we have $\lim_{n\to\infty} \|P^{nD+r+D-i}(x,\cdot) - \pi_r(\cdot)\| = 0$. The result follows by averaging over $1 \leq r \leq D$ and using the triangle inequality. □

The above conclusions still allow for the possibility of a null set $G^C$ from which convergence fails. This null set can indeed arise, even for simple examples, on both discrete and continuous state spaces:

EXAMPLE 3 ([5, 17]). Let $\mathcal{X} = \{1, 2, \ldots\}$. Let $P(1, \{1\}) = 1$, and for $x \geq 2$, $P(x, \{1\}) = 1/x^2$ and $P(x, \{x+1\}) = 1 - (1/x^2)$. This chain has stationary distribution $\pi(\cdot) = \delta_1(\cdot)$ and it is $\phi$-irreducible (with respect to $\pi$) and aperiodic. On the other hand, if $X_0 = x \geq 2$, then $\mathbf{P}[X_n = x+n$ for all $n] = \prod_{j=x}^{\infty}(1 - (1/j^2)) > 0$ so that $\|P^n(x,\cdot) - \pi(\cdot)\| \not\to 0$. Hence, convergence holds only from the set $G = \{1\}$, but fails to hold from the set $G^C = \{2, 3, 4, \ldots\}$. Of course, this example is not irreducible in the classical sense since no state $x \geq 2$ is reachable from the state 1. However, it is still *indecomposable* (see, e.g., [21]).

EXAMPLE 4 (Continuous state space version). Let $\mathcal{X} = [0,1]$. Define the transition kernel $P(x, \cdot)$ as follows. If $x = 1/m$ for some positive integer $m$, then $P(x, \cdot) = x^2 \operatorname{Uniform}[0, 1] + (1-x^2)\delta_{1/(m+1)}(\cdot)$. For all other $x$, $P(x, \cdot) = \operatorname{Uniform}[0, 1]$. Then the chain has stationary distribution $\pi(\cdot) = \operatorname{Uniform}[0,1]$ and it is $\phi$-irreducible (with respect to $\pi$) and aperiodic. On the other hand, if $X_0 = 1/m$ for some $m \geq 2$, then $\mathbf{P}[X_n = 1/(m+n)$ for all $n] = \prod_{j=m}^{\infty}(1-(1/j^2)) > 0$ so that $\|P^n(x,\cdot) - \pi(\cdot)\| \not\to 0$. Hence, convergence fails to hold from the set $G^C = \{1/2, 1/3, 1/4, \ldots\}$.

To rectify the problems of the null set $G^C$, we consider Harris recurrence, a concept developed by Harris [8] and introduced to statisticians in the important works of Tierney [25] and Chan and Geyer [3] (see also [6]).

DEFINITION 5. A $\phi$-irreducible Markov chain with stationary distribution $\pi(\cdot)$ is *Harris recurrent* if for all $A \subseteq \mathcal{X}$ with $\pi(A) > 0$ and all $x \in \mathcal{X}$, we have $\mathbf{P}(\tau_A < \infty | X_0 = x) = 1$.

We now prove a number of equivalences.

THEOREM 6. *For a $\phi$-irreducible Markov chain with stationary probability distribution $\pi(\cdot)$ and period $D \geq 1$, the following are equivalent:*

(i) *The chain is Harris recurrent.*
(ii) *For all $A \subseteq \mathcal{X}$ with $\pi(A) > 0$ and all $x \in \mathcal{X}$, we have $\mathbf{P}(X_n \in A \text{ i.o.} | X_0 = x) = 1$. (Here i.o. means "infinitely often," i.e., for infinitely many different times $n$.)*
(iii) *For all $x \in \mathcal{X}$, $\lim_{n \to \infty} \|(1/D) \sum_{r=1}^{D} P^{nD+r}(x, \cdot) - \pi(\cdot)\| = 0$.*
(iv) *For all $x \in \mathcal{X}$, $\mathbf{P}[\tau_G < \infty | X_0 = x] = 1$, where $G$ is as in Proposition 2.*
(v) *For all $x \in \mathcal{X}$ and all $A \in \mathcal{F}$ with $\pi(A) = 1$, $\mathbf{P}[\tau_A < \infty | X_0 = x] = 1$.*
(vi) *For all $x \in \mathcal{X}$ and all $A \in \mathcal{F}$ with $\pi(A) = 0$, $\mathbf{P}[X_n \in A \text{ for all } n | X_0 = x] = 0$.*

PROOF. (ii) $\Longrightarrow$ (i); (i) $\Longrightarrow$ (v); (v) $\Longrightarrow$ (iv); and (v) $\Longleftrightarrow$ (vi): Immediate.

(i) $\Longrightarrow$ (ii): Suppose to the contrary that (ii) does not hold. Then there is some $A \subseteq \mathcal{X}$ with $\pi(A) > 0$, some $x \in \mathcal{X}$ and some $N \in \mathbf{N}$ such that $\mathbf{P}(X_n \notin A \ \forall n \geq N | X_0 = x) > 0$. Integrating over choices of $y = X_N$, this implies there is some $y \in \mathcal{X}$ with $\mathbf{P}(\tau_A = \infty | X_0 = y) > 0$, contradicting (i).

(iv) $\Longrightarrow$ (iii): From Proposition 2, once the chain reaches $G$, it will converge. The convergence in (iii) then follows. More formally, conditional on the first hitting time $\tau_G$ and the corresponding chain value $X_{\tau_G}$, the chain will converge in total variation distance as in (iii). Statement (iii) then follows by integrating over all choices of $\tau_G$ and $X_{\tau_G}$ and using the triangle inequality for total variation distance.

(iii) $\Longrightarrow$ (i): If $\phi(A) > 0$ [where $\phi(\cdot)$ is an irreducibility measure], then we must have $\pi(A) > 0$ (see, e.g., Lemma 3 of [23]), so by (iii) we have that for all $x \in \mathcal{X}$, $\sum_{r=1}^{D} P^{nD+r}(x, A) \to D\pi(A) > 0$ and, in particular, $\sum_{n=1}^{\infty} P^n(x, A) = \infty$. It then follows from Theorem 9.0.1 of [13], using their definition of recurrence on page 182, that we can find an absorbing subset $H \subseteq \mathcal{X}$ such that the chain restricted to $H$ is Harris recurrent and $\pi(H) = 1$. Then $(1/D) \sum_{r=1}^{D} P^{nD+r}(x, H) \to 1$, so $P^n(x, H) \to 1$. Hence, the chain will eventually reach $H$ with probability 1. Since the chain restricted to $H$ is Harris recurrent, it will then eventually reach any $A$ with $\pi(A) = 1$, with probability 1, thus establishing (i). □

For completeness, we note another method for verifying Harris recurrence (although we do not use it here). Given a Markov chain with stationary distribution $\pi(\cdot)$, a subset $C \in \mathcal{F}$ is *small* if $\pi(C) > 0$ and there is an $\varepsilon > 0$ and a probability measure $\nu(\cdot)$ on $(\mathcal{X}, \mathcal{F})$ such that $P(x, A) \geq \varepsilon \nu(A)$ for all $A \in \mathcal{F}$ and $x \in C$. It easily follows that we must have $\nu \ll \pi$.

PROPOSITION 7 ([13]). *If a Markov chain with stationary distribution $\pi(\cdot)$ has a small set $C$ with the property that $\mathbf{P}[X_n \in C \text{ i.o.} | X_0 = x] = 1$ for all $x \in \mathcal{X}$, then the chain is Harris recurrent.*

PROOF. Using the splitting technique (see, e.g., [15]), each time it is in $C$, we can regard the chain as proceeding by moving according to $\nu(\cdot)$ with probability $\varepsilon$. If the chain returns to $C$ infinitely often, then with probability 1 it will eventually move according to $\nu(\cdot)$. Since $\nu \ll \pi$, this means it will eventually leave any set of null $\pi$-measure. Hence, the result follows from Theorem 6(vi). □

Various drift conditions can be used to establish that $\mathbf{P}[X_n \in C \text{ i.o.} | X_0 = x] = 1$ for all $x \in \mathcal{X}$ and thus establish Harris recurrence. For example, it follows from [13], Theorem 13.0.1, that for $\phi$-irreducible chains, it suffices that there exists a measurable function $V : \mathcal{X} \to (0, \infty)$ such that $\mathbf{E}[V(X_1) | X_0 = x] \leq V(x) - 1 + b\mathbf{1}_C(x)$ for some $b < \infty$. Alternatively, it follows from [13], Theorem 8.4.3 (see also [6]) that for $\phi$-irreducible chains, Harris recurrence follows if there exists a measurable function $V : \mathcal{X} \to (0, \infty)$ such that $V^{-1}((0, \alpha])$ is small for all $\alpha > 0$ and such that $\mathbf{E}[V(X_1) | X_0 = x] \leq V(x)$ for all $x \in \mathcal{X} \setminus C$.

REMARK. We note that the null sets related to Harris recurrence are of an "extreme" kind in the sense that the chain may fail entirely to converge from the null set. Less radically, one could consider chains which converge from everywhere but which have a slower qualitative *rate* of convergence from some null set. For example, it should be possible to construct Markov chains which are Harris recurrent and geometrically ergodic but which converge at a subgeometric rate from a certain null set of initial states; or, chains which are polynomially ergodic at a particular polynomial rate $\alpha$ but which fail to converge at the polynomial rate $\alpha$ from some null set; or, chains which are geometrically ergodic but which fail to converge from one null set, converge polynomially from another null set, converge subpolynomially from a third null set, and so on. In this context, Harris recurrence can be seen as one in a series of properties ensuring that "things are not worse when starting from a null set than when starting from anywhere else."

**3. Full-dimensional Metropolis–Hastings algorithms.** Let $\mathcal{X}$ be some state space with $\sigma$-algebra $\mathcal{F}$. Let $\pi(\cdot)$ be a probability distribution on $(\mathcal{X}, \mathcal{F})$ having unnormalized density function $f : \mathcal{X} \to (0, \infty)$ with respect to some reference measure $\nu(\cdot)$ so that $\int_\mathcal{X} f(x)\nu(dx) < \infty$ and

$$\pi(A) = \frac{\int_A f(x)\nu(dx)}{\int_\mathcal{X} f(x)\nu(dx)}, \qquad A \in \mathcal{F}.$$

Note that we assume $f > 0$ on $\mathcal{X}$ or, equivalently, that $\mathcal{X}$ is defined to be the support of $f$. To avoid trivialities, we assume that $\pi(\cdot)$ is not concentrated at a single state, that is, that $\pi\{x\} < 1$ for all $x \in \mathcal{X}$.

Let $q : \mathcal{X} \times \mathcal{X} \to [0, \infty)$ be any jointly measurable function such that $\int q(x, y)\nu(dy) = 1$ for all $x \in \mathcal{X}$. Define the Markov kernel $Q(x, \cdot)$ by $Q(x, A) = \int_A q(x, y)\nu(dy)$ for $x \in \mathcal{X}$ and let

$$\alpha(x, y) = \min\left[1, \frac{f(y)q(y, x)}{f(x)q(x, y)}\right], \qquad x, y \in \mathcal{X}$$

[with $\alpha(x, y) = 1$ if $f(x)q(x, y) = 0$].

The *full-dimensional Metropolis–Hastings algorithm* [9, 12, 25] proceeds as follows. Given that the chain is in state $X_n$ at time $n$, it generates a "proposal state" $Y_{n+1} \sim Q(X_n, \cdot)$. Then, with probability $\alpha(X_n, Y_{n+1})$, it "accepts" this proposal and sets $X_{n+1} = Y_{n+1}$; otherwise, with probability $1 - \alpha(X_n, Y_{n+1})$, it "rejects" this proposal and sets $X_{n+1} = X_n$. It is easy to check that $\pi(\cdot)$ is then stationary for the Markov chain $\{X_n\}$.

Clearly, any such Markov chain can be decomposed as

$$P(x, A) = (1 - r(x))M(x, A) + r(x)\delta_x(A), \qquad x \in \mathcal{X}, \ A \subseteq \mathcal{X},$$

where $\delta_x(\cdot)$ is a point-mass at $x$, $r(x) = \int q(x, y)[1 - \alpha(x, y)]\nu(dy)$ is the probability of rejection when starting at $X_n = x$ and $M(x, \cdot)$ is the kernel conditional on moving (i.e., on $X_{n+1} \neq X_n$). In particular, the probability distribution $M(x, \cdot)$ is absolutely continuous with respect to $\nu(\cdot)$ for all $x \in \mathcal{X}$.

Regarding Harris recurrence, we have the following result, which was originally proved by Tierney [25] using the theory of harmonic functions:

THEOREM 8 ([25]). *Every $\phi$-irreducible, full-dimensional Metropolis–Hastings algorithm is Harris recurrent.*

PROOF. Since the chain is $\phi$-irreducible and $\pi\{x\} < 1$, we must have $r(x) < 1$ for all $x \in \mathcal{X}$. Suppose $\pi(A) = 1$. Then $\pi(A^C) = 0$, and so, as we are assuming that $f > 0$ throughout $\mathcal{X}$, we also have $\nu(A^C) = 0$. Hence, by absolute continuity, $M(x, A^C) = 0$, that is, $M(x, A) = 1$. It follows that, if the chain is at $x$, then it will eventually move according to $M(x, \cdot)$, at which point it will necessarily move into $A$. The result then follows from Theorem 6(v). □

**4. Metropolis-within-Gibbs algorithms.** We now define Metropolis-within-Gibbs Markov chains [12].

For simplicity, let $\mathcal{X}$ be an open subset of $\mathbf{R}^d$ with Borel $\sigma$-algebra $\mathcal{F}$ and (unnormalized) target density $f : \mathcal{X} \to (0, \infty)$ with $\int_\mathcal{X} f(x) \lambda(dx) < \infty$ [where $\lambda(\cdot)$ is $d$-dimensional Lebesgue measure]. For $1 \le i \le d$, let $q_i : \mathcal{X} \times \mathbf{R} \to [0, \infty)$ be jointly measurable with $\int_{-\infty}^\infty q_i(x, z)\, dz = 1$ for all $x \in \mathcal{X}$ (where $dz$ is one-dimensional Lebesgue measure).

Let $Q_i(x, \cdot)$ be the Markov kernel on $\mathbf{R}^d$ which replaces the $i$th coordinate by a draw from the density $q_i(x, \cdot)$, but leaves the other coordinates unchanged. That is,

$$Q_i(x, S_{i,a,b}) = \int_a^b q_i(x, z)\, dz,$$

where

$$S_{i,a,b} = \{y \in \mathcal{X} : y_j = x_j \text{ for } j \ne i \text{ and } a \le y_i \le b\}.$$

To avoid technicalities and special cases, assume that $Q_i(x, \mathcal{X}) > 0$ for all $x \in \mathcal{X}$ and $1 \le i \le d$, and also that each $q_i$ is *symmetrically positive* in the sense that

$$q_i((x_1, \ldots, x_{i-1}, y, x_{i+1}, \ldots, x_d), z) > 0$$
$$\iff q_i((x_1, \ldots, x_{i-1}, z, x_{i+1}, \ldots, x_d), y) > 0.$$

For $x, y \in \mathbf{R}^d$ and $1 \le i \le d$, let

$$\alpha_i(x, y) = \min\left[1, \frac{f(y) q_i(y, x)}{f(x) q_i(x, y)}\right]$$

[with $\alpha_i(x, y) = 1$ if $f(x) q_i(x, y) = 0$]. Let $P_i$ be the kernel which proceeds as follows. Given $X_n$, it generates a proposal $Y_{n+1} \sim Q_i(X_n, \cdot)$. Then, with probability $\alpha_i(X_n, Y_{n+1})$, it accepts this proposal and sets $X_{n+1} = Y_{n+1}$; otherwise, with probability $1 - \alpha_i(X_n, Y_{n+1})$, it rejects the proposal and sets $X_{n+1} = X_n$.

In terms of these definitions, the Metropolis-within-Gibbs Markov chain proceeds as follows. Random variables $I_1, I_2, \ldots$ taking values in $\{1, 2, \ldots, d\}$ are chosen according to some scheme. (The two most common schemes are random-scan, where $\{I_n\}$ are i.i.d. uniform on $\{1, 2, \ldots, d\}$, and deterministic-scan, where $I_1 = 1, I_2 = 2, \ldots, I_d = d, I_{d+1} = 1, \ldots$.) Then for $n = 0, 1, 2, \ldots$, given $X_n$, the chain generates $X_{n+1} \sim P_{I_{n+1}}(X_n, \cdot)$. It is straightforward to verify that this chain has stationary distribution $\pi(\cdot)$ given by

$$\pi(A) = \frac{\int_A f(x) \lambda(dx)}{\int_\mathcal{X} f(x) \lambda(dx)}, \qquad A \in \mathcal{F}.$$

The above description defines Metropolis-within-Gibbs chains as we shall study them. We can now ask, under what conditions are such chains Harris recurrent? One might think that a result similar to Theorem 8 holds for Metropolis-within-Gibbs algorithms, at least when the target and proposal densities are continuous. However, surprisingly, this is false:

EXAMPLE 9. We present a Metropolis-within-Gibbs algorithm on an open subset $\mathcal{X} \subseteq \mathbf{R}^2$, with stationary distribution $\pi(\cdot)$, with continuous target and proposal densities, which is $\phi$-irreducible and aperiodic, but which fails to converge in distribution to $\pi(\cdot)$ from an uncountable number of different starting points (having total $\pi$-measure zero, of course).

Let $\mathcal{X} = \{(x_1, x_2) \in \mathbf{R}^2 : x_1 > 1\}$ and define the function $f : \mathcal{X} \to (0, \infty)$ by $f(x_1, x_2) = (e/2) \exp(x_1 - |x_2|e^{2x_1})$ (so that $\int_{\mathcal{X}} f = 1$). Let $Q_1(x, \cdot)$ and $Q_2(x, \cdot)$ be symmetric unit normal proposals so that $q_i(x, z) = (2\pi)^{-1/2} \exp(-(z - x_i)^2/2)$ ($i = 1, 2$). Then, clearly, $f$, $q_1$ and $q_2$ are positive continuous functions; it follows that the chain is $\phi$-irreducible where $\psi$ is Lebesgue measure.

Consider the random-scan (say) Metropolis-within-Gibbs Markov chain corresponding to these choices. We shall prove that this chain is *not* Harris recurrent. Indeed, let $S = \{(x_1, 0) : x_1 > 1\}$ be the part of the line $\{x_2 = 0\}$ which lies in $\mathcal{X}$. We claim that if the chain starts at any initial state in $S$, then there is positive probability that it will drift off to $x_1 \to \infty$ without ever updating $x_2$, that is, without ever leaving $S$. Then, since $\pi(S) = 0$, it follows that if $X_0 \in S$, the chain will fail to converge to $\pi(\cdot)$.

To establish the claim, consider first a Markov chain $\{W_n\}$ equivalent to just the first coordinate of $X_n$, under just the kernel $P_1$ (which proposes moves only in the $x_1$ direction), restricted to the state space $S$. Now, on $S$, the density $f$ is proportional to $e^{x_1}$. It follows that for any $\delta > 0$ and $x_1 > 1$, $\alpha_1((x_1, 0), (x_1 - \delta, 0)) \leq e^{-\delta}$, while $\alpha_1((x_1, 0), (x_1 + \delta, 0)) = 1$. That is, proposals to increase $x_1$ will all be accepted, while a positive fraction of the proposals to decrease $x_1$ will be rejected. It follows that on $S$, the kernel $P_1$ has positive drift. Hence, there exists $c > 0$ such that for all $x_1 > 1$,

$$\mathbf{P}[W_n \geq cn \text{ for all sufficiently large } n | W_0 = x_1] > 0.$$

On the other hand, the density $f(x_1, x_2)$ as a function of $x_2$ alone (i.e., with $x_1$ regarded as a constant) is proportional to $\exp(-|x_2|e^{2x_1})$. It follows that the probability of accepting a proposal in the $x_2$ direction is equal to $\mathbf{E}[\exp(-|Z|e^{2x_1})]$, where $Z \sim N(0, 1)$, which is less than

$$\int_{-\infty}^{\infty} \exp(-|z|e^{2x_1}) \, dz = 2e^{-2x_1}.$$

Now, since $\sum_n 2e^{-2cn} < \infty$, it follows from the Borel–Cantelli lemma (e.g., [22], Theorem 3.4.2) that there is positive probability that all proposed moves in the $x_2$ direction will be rejected. That is, for any $x_1 > 1$,

$$\mathbf{P}[X_n \in S \text{ for all } n | X_0 = (x_1, 0)] > 0,$$

thus proving the claim. (We shall see in Corollary 18 below that the "problem" with this example is that the one-dimensional integral of $f$ over the line $\{x_2 = 0\}$ is infinite.)

REMARK. In the above example, it is also possible, if desired, to modify $f$ to decrease to 0 near the boundary $\{x_1 = 1\}$ in order to make $f$ be continuous throughout $\mathbf{R}^2$ (not just on $\mathcal{X}$) without affecting the result.

To proceed, decompose $P_i$ as $P_i(x, \cdot) = [1 - r(x)]M_i(x, \cdot) + r(x)\delta_x(\cdot)$ so that $M_i(x, S)$ is the kernel corresponding to moving (i.e., both proposing and accepting) in the $i$th direction.

LEMMA 10. *Let $(i_1, i_2, \ldots, i_n)$ be any sequence of coordinate directions. Assume that each of the $d$ coordinate directions appears at least once in the sequence $(i_1, i_2, \ldots, i_n)$. Then $M_{i_1} M_{i_2} \cdots M_{i_n}$ is an absolutely continuous kernel, that is, if $A \in \mathcal{F}$ with $\lambda(A) = 0$, then $(M_{i_1} M_{i_2} \cdots M_{i_n})(x, A) = 0$ for all $x \in \mathcal{X}$.*

PROOF. We shall compute a density for $(M_{i_1} M_{i_2} \cdots M_{i_n})(x, \cdot)$. The result then follows since every distribution having a density is absolutely continuous. Let

$$J = \{m : 1 \leq m \leq n, i_j \neq i_m \text{ for } m < j \leq n\},$$

that is, $J$ is the set of "last time the chain moved in direction $i$" for each coordinate $i$. (Thus, $|J| = d$.) For $1 \leq m \leq n$, let $S_m \subseteq \mathbf{R}$ be any Borel subset so that $S = S_1 \times \cdots \times S_d$ is an arbitrary measurable rectangle in $\mathbf{R}^d$. Then define subsets $R_m \subseteq \mathbf{R}$ for $1 \leq m \leq n$ by letting $R_m = S_{i_m}$ if $m \in J$, and $R_m = \mathbf{R}$ otherwise.

We then compute that

$$(M_{i_1} M_{i_2} \cdots M_{i_n})(x, S)$$
$$= \int_{R_1} \int_{R_2} \cdots \int_{R_n} q_{i_1}(x_1, x_2)\alpha(x_1, x_2) q_{i_2}(x_2, x_3)\alpha(x_2, x_3) \times \cdots$$
$$\times q_{i_{n-1}}(x_{n-1}, x_n)\alpha(x_{n-1}, x_n) \, dx_1 \, dx_2 \cdots dx_n.$$

It follows that the density of $(M_{i_1} M_{i_2} \cdots M_{i_n})(x, \cdot)$ is given by the above formula, but with the integration over the variables $\{x_j; j \in J\}$ omitted. Hence, $(M_{i_1} M_{i_2} \cdots M_{i_n})(x, \cdot)$ has a density and is thus absolutely continuous. □

From the law of total probability, we therefore obtain the following:

COROLLARY 11. *If $A$ has Lebesgue measure 0, then $\mathbf{P}[X_n \in A | X_0 = x_0] \leq \mathbf{P}[D_n]$ where $D_n$ is the event that by time $n$, the chain has not yet moved in each coordinate direction.*

This allows us to prove the following:

THEOREM 12. *Consider a $\phi$-irreducible Metropolis-within-Gibbs Markov chain. Suppose that from any initial state $x$, with probability 1, the chain will eventually move at least once in each coordinate direction. Then the chain is Harris recurrent.*

PROOF. The hypothesis implies that for all $x \in \mathcal{X}$, $\lim_{n \to \infty} \mathbf{P}[D_n | X_0 = x] = 0$. Now, let $\pi(A) = 0$. Then since $f > 0$ on $\mathcal{X}$, we must also have $\lambda(A) = 0$. Hence, by Corollary 11, we must have

$$\mathbf{P}[X_n \in A \ \forall n | X_0 = x] \leq \lim_{n \to \infty} \mathbf{P}[X_n \in A | X_0 = x] \leq \lim_{n \to \infty} \mathbf{P}[D_n | X_0 = x] = 0.$$

The result then follows from Theorem 6(vi). □

The classical Gibbs sampler (see, e.g., [4]) is a special case of Metropolis-within-Gibbs in which the proposal densities are chosen so that $\alpha(x, y) \equiv 1$, that is, so that all proposed moves are accepted. Now, with either the deterministic-scan or systematic-scan Gibbs sampler variants, it is certainly true that with probability 1, moves are eventually proposed in all directions. So, since $\alpha(x, y) \equiv 1$, it also follows that with probability 1, the chain will eventually move in all directions. Hence, from Theorem 12, we immediately obtain the following:

COROLLARY 13 ([25]). *Every $\phi$-irreducible deterministic- or random-scan Gibbs sampler is Harris recurrent.*

**5. Subchains of Metropolis-within-Gibbs algorithms.** We now consider the extent to which Harris recurrence of the full chain can be "inherited" from Harris recurrence of various subchains. For a subset $I = \{i_1, \ldots, i_r\} \subseteq \{1, \ldots, n\}$, let $P^{[I]}$ be the Markov kernel which corresponds to the original Metropolis-within-Gibbs chain, that it is except conditional on never moving in any coordinate directions other than the coordinate directions $i_1, \ldots, i_r$. Call the collection of kernels $P^{[I]}$, where $|I| = d - 1$, the "$(d-1)$-dimensional subchains." These subchains can fail to be $\phi$-irreducible:

EXAMPLE 14. Suppose that

$$\mathcal{X} = \{(x_1, x_2) \in \mathbf{R}^2 : 16 < x_1^2 + x_2^2 < 25\}$$

(an annulus or "donut-shaped" state space) and that the proposal kernels $Q_i(x,\cdot)$ simply replace $x_i$ by a draw from the Uniform$[x_i - 1, x_i + 1]$ distribution. Then the full Metropolis-within-Gibbs chain is $\phi$-irreducible, but the one-dimensional subchain along the line $\{x_2 = 0\}$ breaks up into two distinct noncommunicating intervals, $(-5,-4)$ and $(4,5)$, and is therefore not $\phi$-irreducible.

Harris recurrence is often defined solely for $\phi$-irreducible chains (e.g., [13]). We generalize as follows. Call a chain *piecewise Harris* if the state space $\mathcal{X}$ can be partitioned into a disjoint union $\mathcal{X} = \bigcup_{\alpha \in \mathcal{S}} \mathcal{X}_\alpha$ where each $\mathcal{X}_\alpha$ is closed and the chain restricted to each $\mathcal{X}_\alpha$ is Harris recurrent. Of course, if the partition consists solely of a single $\mathcal{X}_\alpha$, then the full chain is Harris recurrent. The following proposition says that the piecewise Harris property often suffices:

PROPOSITION 15. *If a Markov chain is piecewise Harris and is also $\phi$-irreducible, then it is Harris recurrent.*

PROOF. Let $\mathcal{X}_\beta$ be any nonempty element of the partition from the definition of piecewise Harris and let $\mathcal{X}_* = \mathcal{X} \setminus \mathcal{X}_\beta$. If $\phi(\mathcal{X}_*) > 0$, then, by $\phi$-irreducibility, for each $x \in \mathcal{X}$, there exists $n = n(x)$ with $P^n(x, \mathcal{X}_*) > 0$. Since $\mathcal{X}_\beta$ is closed, this implies that $\mathcal{X}_\beta$ is empty, contradicting our assumption. Hence, $\phi(\mathcal{X}_*) = 0$. Since $\phi$ is nonzero, we must have $\phi(\mathcal{X}_\beta) > 0$. It then follows similarly that $\mathcal{X}_*$ is empty, that is, that $\mathcal{X}_\beta = \mathcal{X}$. Thus, the partition contains just a single element, and so the chain is Harris recurrent. $\square$

In terms of the piecewise Harris property, we have the following:

THEOREM 16. *Consider a random-scan Metropolis-within-Gibbs chain, as above. Suppose that all the $(d-1)$-dimensional subchains in every $(d-1)$-dimensional coordinate hyperplane are piecewise Harris. Then the full chain is piecewise Harris. (In particular, by Proposition 15, if the full chain is $\phi$-irreducible, then the full chain is Harris recurrent.)*

PROOF. Consider any fixed initial state $x_0 = (x_{0,1}, \ldots, x_{0,d})$. By Theorem 12, it suffices to show that, with probability 1, when starting at $X_0 = x_0$, the chain will eventually move in each coordinate direction.

Suppose, to the contrary, that this is false and that there is positive probability that the chain never moves in some direction, say (for notational simplicity) in direction $d$. Let $H = \{y \in \mathcal{X} : y_j = x_{0,j} \text{ for } j \neq d\}$ be the hyperplane corresponding to never moving in the $d$th direction.

Let $I_n$ be the direction of the proposed move of the full chain at time $n$ and let $A_n = 1$ if this move is accepted, otherwise let $A_n = 0$. Let

$$C_{m,r} = \{w \in H; \mathbf{P}[I_m = d, A_m = 1 | X_0 = w] \geq 1/r\}.$$

That is, $C_{m,r}$ is the set of states in $H$ which have probability $\geq 1/r$ of changing the $d$th coordinate, $m$ steps later, when moving according to the subchain.

By assumption, $Q_d(x, \mathcal{X}) > 0$ and $f(x) > 0$ for all $x \in \mathcal{X}$. This implies that the chain has positive probability, starting from any $x \in H$, of eventually moving in the $d$th direction, that is, of leaving $H$. Hence,

$$\bigcup_{m,r=1}^{\infty} C_{m,r} = H. \tag{1}$$

(In fact, it suffices to consider just $m = 1$, but we do not use that fact here.)

Consider now the subchain $P^{[1,2,\ldots,d-1]}$, restricted to the hyperplane $H$. Since this subchain is piecewise Harris, we must have $x_0 \in H_0$ for some closed subset $H_0 \subseteq H$ such that the subchain restricted to $H_0$ is Harris recurrent with respect to some nonzero measure $\psi_j(\cdot)$. From (1), there must exist some $m, r \in \mathbf{N}$ with $\psi_j(C_{m,r}) > 0$. It then follows from Theorem 6(ii) that, with probability 1, $C_{m,r}$ is hit infinitely often by the subchain. In other words, conditional on the full chain never moving in the $d$th direction, it will enter $C_{m,r}$ infinitely often.

However, each time the full chain visits $C_{m,r}$, it has probability $\geq 1/r$ of moving in the $d$th direction $m$ steps later. It follows that, with probability 1, the full chain will eventually jump in the $d$th direction and hence leave $H$. This contradicts our assumption that the chain has positive probability of never leaving $H$. □

Unfortunately, Theorem 16 still requires that we verify Harris recurrence of various subchains, which may be difficult. However, if the subchains of all dimensions all have stationary distributions, then no Harris recurrence needs to be checked as we see in the following:

COROLLARY 17. *Consider a random-scan Metropolis-within-Gibbs chain as above. Suppose that every $r$-dimensional subchain in every $r$-dimensional coordinate hyperplane has a each have stationary probability distribution, for all $1 \leq r \leq d$. Then the full chain is piecewise Harris (as are all the subchains).*

PROOF. Let $T_r$ be the statement that all the subchains of dimension $\leq r$ are piecewise Harris. Then $T_1$ holds by Theorem 8. Furthermore, from Theorem 16, for any $r < d$, if $T_r$ holds, then $T_{r+1}$ must hold. Hence, the result follows by induction. □

We then have the following:

COROLLARY 18. *Consider a random-scan Metropolis-within-Gibbs Markov chain. Suppose that the target density $f$ has the property that its $r$-dimensional integral has finite Lebesgue integral, over every $r$-dimensional coordinate hyperplane of $\mathcal{X}$, for all $1 \le r \le d$. Then the full chain is piecewise Harris (as are all the subchains).*

PROOF. In this case, $f$ is an (unnormalized) density for a stationary probability distribution of each subchain on each hyperplane. (Note that the Lebesgue integral of $f$ over the hyperplane must be positive since we are assuming that $\mathcal{X}$ is open and that $f > 0$ on $\mathcal{X}$.) Hence, the result follows from Corollary 17. □

In the counterexample of Example 9, the one-dimensional $x_1$-chain fails to have a stationary distribution along the line $\{x_2 = 0\}$ since the integral of $f$ along the line $\{x_2 = 0\}$ is infinite.

A result similar to Corollary 18 appears as Theorem 1 of [3] under the assumption that each subchain is $\phi$-irreducible (which, as we have seen in Example 14, can easily fail to hold):

COROLLARY 19 ([3]). *Consider a random-scan Metropolis-within-Gibbs Markov chain. Suppose that the target density $f$ has the property that its $r$-dimensional integral has finite Lebesgue integral over every $r$-dimensional coordinate hyperplane of $\mathcal{X}$, for all $1 \le r \le d$. If the full chain and all the subchains are $\phi$-irreducible, then the full chain is Harris recurrent.*

**6. Trans-dimensional MCMC algorithms.** In certain statistical setups (e.g., autoregressive models), the number of parameters is not fixed in advance. This means that the state space of possible parameter values is a (disjoint) union of spaces of different dimensions. Exploring such state spaces through MCMC algorithms requires the introduction of trans-dimensional MCMC. Trans-dimensional MCMC algorithms first appeared in [14] and [16]; their introduction into modern statistical practice (under the name "reversible jump") is due to the influential paper of Green [7] (see also [26]).

Suppose that for each $m \in \mathcal{M} \subseteq \mathbf{N}$, where $|\mathcal{M}| > 1$ (and usually $|\mathcal{M}| = \infty$), we have a space $\mathcal{X}_m$ of dimension $d_m$, that is, $\mathcal{X}_m$ is an open subset of $\mathbf{R}^{d_m}$. We combine these different spaces into a single state space $\mathcal{X}$ by setting $\mathcal{X} = \bigcup_{m=1}^{\infty} (\{m\} \times \mathcal{X}_m)$. Furthermore, suppose that on each $\mathcal{X}_m$, we have an unnormalized target density function $f_m : \mathcal{X}_m \to (0, \infty)$ with $\int_{\mathcal{X}_m} f_m < \infty$. We then combine that into a single probability distribution $\pi(\cdot)$ on $\mathcal{X}$ by choosing some $p : \mathcal{M} \to (0,1)$ with $\sum_{m \in \mathcal{M}} p(m) = 1$, setting

$$\pi(m, A) = p(m) \frac{\int_A f_m(x) \lambda_m(dx)}{\int_{\mathcal{X}} f_m(x) \lambda_m(dx)} \tag{2}$$

and using linearity; in (2), $\lambda_m(\cdot)$ is Lebesgue measure on $\mathbf{R}^{d_m}$.

Trans-dimensional chains can proceed in a variety of ways [2, 7]. We first consider a general class which we call *full-dimensional trans-dimensional MCMC*. Fix some $0 < a < 1$ and some irreducible kernel $R(m, \cdot)$ on $\mathcal{M}$ such that $R(m, m') > 0$ if and only if $R(m', m) > 0$. Then at each iteration, with probability $a$, the chain proposes a "between-model move" which replaces $(m, x)$ by $(m', x')$, where $m' \sim R(m, \cdot)$ and where $x' \in \mathcal{X}_{m'}$ is generated by some complicated dimension-matching scheme [7]. This proposal is then accepted or rejected according to the usual Metropolis–Hastings scheme, except that now the formula for $\alpha[(m, x), (m', x')]$ is more complicated and involves a Jacobian of the transformation used to generate $x'$. Otherwise, with probability $1 - a$, the chain leaves $m$ fixed but proposes a "within-model move," that is, to replace $x$ by $x' \in \mathcal{X}_m$, using a full-dimensional Metropolis–Hastings proposal on $\mathcal{X}_m$.

What about Harris recurrence? We have the following:

THEOREM 20. *Consider a full-dimensional trans-dimensional MCMC algorithm as above. Let $D$ be the event that no within-model move is ever accepted. Suppose that $\mathbf{P}[D|X_0 = (m, x)] = 0$ for all $(m, x) \in \mathcal{X}$. Then the algorithm is Harris recurrent.*

PROOF. The proof is very similar to that of Theorem 8. Since $\mathbf{P}[D|X_0 = (m, x)] = 0$, the chain must eventually accept a within-model move. But since the within-model proposal distributions are full-dimensional, the probability of remaining in any set of $\pi$-measure 0 after such a move is equal to 0. The result thus follows from Theorem 6(vi). □

REMARK. Theorem 20 remains true regardless of the details of how the between-model moves are implemented, provided only that they preserve the stationarity of $\pi(\cdot)$.

REMARK. Even without verifying the hypothesis of Theorem 20, it is true that once a full-dimensional trans-dimensional MCMC algorithm makes at least one within-model move, then, since the within-model moves are full-dimensional, with probability 1, the chain will move to the set $G$ of Proposition 2 and hence will then converge. The issue in Theorem 20 is whether or not such a within-model move will eventually occur with probability 1.

Now, Theorem 20 allows for the possibility that from a null set, the model numbers $m_n$ might have positive probability of converging to $+\infty$ without ever accepting a within-model move. This seems quite plausible. On the other hand, if the $\{m_n\}$ process is recurrent, the situation is less clear due to the complicated details of the $(m, x) \to (m', x')$ mapping corresponding to the between-model moves. Conditional on never accepting a within-model move, even if the chain returns to $\mathcal{X}_1$ (say) infinitely often, it might potentially be at "worse and worse" points within $\mathcal{X}_1$ each time it returned and thus have an ever decreasing probability of accepting within-model moves. So, even a chain in which $\{m_n\}$ is recurrent could conceivably fail to be Harris recurrent. We state this as an open problem:

OPEN PROBLEM 1. Does there exist a $\phi$-irreducible full-dimensional trans-dimensional MCMC algorithm as above, for which $\mathbf{P}[m_n = 1 \text{ i.o.} | X_0 = (m, x)] = 1$ for all $m \in \mathcal{M}$ and $x \in \mathcal{X}_m$, which fails to be Harris recurrent?

More generally, a trans-dimensional MCMC might not be full-dimensional. That is, the within-model moves might themselves be of Metropolis-within-Gibbs form. To model this, we proceed as in [2]. We replace $\mathcal{X}_m$ by $\tilde{\mathcal{X}}_m \equiv \mathcal{X}_m \times [0,1] \times [0,1] \times \cdots$, with stationary distribution $\tilde{\pi}_m = \pi_m \times \text{Uniform}[0,1] \times \text{Uniform}[0,1] \times \cdots$. We then let $h_{ij} : \tilde{\mathcal{X}}_i \to \tilde{\mathcal{X}}_j$ be deterministic functions defined whenever $R(i,j) > 0$, and such that $h_{ji} = (h_{ij})^{-1}$. The between-model moves are specified by requiring that when the algorithm proposes changing $m$ to $m'$, it simultaneously proposes changing $x$ to $h_{mm'}(x)$.

A special case is when each $h_{ij}$ function is simply the identity function, which is plausible if $\mathcal{X}_m = [0,1]^{d_m}$ for each $m \in \mathcal{M}$. More generally, we consider *coordinate-preserving trans-dimensional MCMC* in which each $h_{ij}$ can be decomposed as

$$(3) \qquad h_{ij} = h_{ij}^{(1)} \times h_{ij}^{(2)} \times \cdots,$$

where each $h_{ij}^{(\ell)} : \mathbf{R} \to \mathbf{R}$ and its inverse are differentiable functions acting solely on the $\ell$th coordinate. That is, the between-model moves modify each coordinate separately. Given a current state $X_n = (m, x)$, the algorithm then proceeds as follows. First, it replaces the coordinates $d_m + 1, d_m + 2, \ldots$ by fresh i.i.d. draws from the Uniform$[0,1]$ distribution. (Of course, in practice, we only need to generate such Uniform$[0,1]$ draws when they are required. But from a theoretical perspective, it is simplest to pretend they are updated at each iteration.) Then, with probability $a$, it proposes a between-model move as above, otherwise, with probability $1 - a$, it chooses one of the coordinates $1, 2, \ldots, d_m$ uniformly at random and executes a Metropolis-within-Gibbs move for that coordinate only.

THEOREM 21. *Consider a $\phi$-irreducible, trans-dimensional MCMC chain which is coordinate-preserving as in* (3). *Suppose that for each $(m, x) \in \mathcal{X}$, when the chain starts at $X_0 = (m, x)$, then, with probability $1$, it eventually accepts at least one move in each of the coordinate directions $1, 2, \ldots, d_m$. Then the chain is Harris recurrent.*

PROOF. The proof is analogous to that of Theorem 12. The only difference is that we do not require a move to be accepted in the coordinate directions $d_m + 1, d_m + 2, \ldots$, nor in the direction corresponding to the model index $m$. To justify this, note that since $\mathcal{M}$ is countable with $p(m) > 0$ for all $m \in \mathcal{M}$, every distribution on $\mathcal{M}$ is absolutely continuous. Also, when starting from $X_0 = (m, x)$, coordinates $d_m + 1, d_m + 2, \ldots$ are drawn from an absolutely continuous (Uniform) distribution automatically. So, in the context of Theorem 12, this is equivalent to having already moved in the coordinate directions $d_m + 1, d_m + 2, \ldots$ and in the direction of $\mathcal{M}$. Then, just as in Theorem 12, the chain will eventually leave any set of zero stationary measure. The result thus follows from Theorem 6(vi). □

This leads to the question of Harris recurrence for trans-dimensional chains which are *coordinate-mixing*, that is, not coordinate-preserving. Unfortunately, this situation is more complicated due to lack of control over the composition of the $h_{ij}$ functions. For manageability, call a trans-dimensional chain *dimension-controlling* if $h_{ij}$ is the identity on coordinate directions $\ell > \max(d_i, d_j)$, that is, if $h_{ij}$ does not mix in any more dimensions than are necessary for dimension-matching. Now, it seems that coordinate-mixing in the $h_{ij}$ should only help the chain to avoid null sets. However, the difficulty is that the between-model moves could, for example, "swap" the values of two coordinates so that updating each coordinate position once could correspond to updating one value twice and the other value not at all. Thus, the situation is unclear and we state this as an open problem:

OPEN PROBLEM 2. Consider a $\phi$-irreducible, coordinate-mixing, dimension-controlling trans-dimensional MCMC chain, as above. Suppose that for each $(m, x) \in \mathcal{X}$, when the chain starts at $X_0 = (m, x)$, then, with probability 1, it eventually accepts at least one move in each of the coordinate directions $1, 2, \ldots, d_m$. Does this imply that the chain is Harris recurrent?

A positive answer to this question would show that general trans-dimensional chains, like other Metropolis-within-Gibbs chains, are Harris recurrent provided they eventually move at least once in each coordinate direction with probability 1.

**Acknowledgments.** We are grateful to Charlie Geyer for piquing our interest in Harris recurrence and to Galin Jones for a helpful conversation. We thank the two anonymous referees for very detailed and insightful reports.

DEPARTMENT OF MATHEMATICS AND STATISTICS
FYLDE COLLEGE
LANCASTER UNIVERSITY
LANCASTER LA1 4YF
ENGLAND
E-MAIL: g.o.robert@lancaster.ac.uk

DEPARTMENT OF STATISTICS
UNIVERSITY OF TORONTO
TORONTO, ONTARIO
CANADA M5S 3G3
E-MAIL: jeff@math.toronto.edu